\documentclass{ws-rmta}

\usepackage{algorithm}
\usepackage{algpseudocode}
\usepackage{graphicx}

\begin{document}

\markboth{Alexander Dubbs and Alan Edelman}
{Infinite Random Matrix Theory, Tridiagonal Bordered Toeplitz Matrices, and the Moment Problem}

%
\catchline{}{}{}{}{}
%

\title{Infinite Random Matrix Theory, Tridiagonal Bordered Toeplitz Matrices, and the Moment Problem}

\author{Alexander Dubbs}

\address{Mathematics, MIT, 77 Massachusetts Avenue\\
Cambridge, MA 02139, United States of America, \\
\email{alex.dubbs@gmail.com}}

\author{Alan Edelman}

\address{Mathematics, MIT, 77 Massachusetts Avenue\\
Cambridge, MA 02139, United States of America, \\
\email{edelman@math.mit.edu}}

\maketitle

\begin{history}
\received{April 7, 2014}
\revised{\today}
\end{history}

\begin{abstract}
The four major asymptotic  level density laws of random matrix theory may all be showcased though their Jacobi parameter representation as having a bordered Toeplitz form. We compare and contrast these laws, completing and exploring their representations in one place. Inspired by the bordered Toeplitz form, we propose an algorithm for the finite moment problem by proposing a solution whose density has a bordered Toeplitz form.
\end{abstract}


\section{Introduction}

Consider the ``big''  laws for asymptotic level densities for various random matrices: 
\begin{center}
\begin{tabular}{lc}  
 Wigner semicircle law & \cite{Wigner1958} \\
 Marchenko-Pastur law &\cite{Marcenko1967} \\
 Kesten-McKay law &\cite{Kesten1959}, \cite{McKay1981} \\
 Wachter law & \cite{Wachter1978} \\  
 \end{tabular}
\end{center}

\begin{table}[h!]
\begin{center}\rule{-.5in}{0in}\begin{tabular}{|lll|} \hline
Measure & Support & Parameters \\ \hline
\begin{tabular}{l}
\rule{-.1 in}{.25in}
Wigner semicircle\\ 
$d\mu_{WS} = \displaystyle\frac{\sqrt{4-x^2}}{2\pi}dx$\\
\end{tabular} & $I_{WS} = [\pm 2]$ & N/A \\
&&\\ \hline
\begin{tabular}{l}
\rule{-.1 in}{.25in}
Marchenko-Pastur\\
$d\mu_{MP} = \displaystyle\frac{\sqrt{(\lambda_+-x)(x-\lambda_-)}}{2\pi x}dx$\\
\end{tabular} & $I_{MP} = [\lambda_-,\lambda_+]$ & $\lambda_\pm = (1\pm\sqrt{\lambda})^2$, $\lambda\geq 1$\\
&&\\ \hline
\begin{tabular}{l}
\rule{-.1 in}{.25in}
Kesten-McKay\\
$d\mu_{KM} = \displaystyle\frac{v\sqrt{4(v-1)-x^2}}{2\pi (v^2-x^2)}dx$\\
\end{tabular} & $I_{M} = [\pm 2\sqrt{v-1}] $ & $v\geq 2$\\
&&\\ \hline
\begin{tabular}{l}
\rule{-.1 in}{.25in}
Wachter\\
$d\mu_W = \displaystyle\frac{(a+b)\sqrt{(\mu_+-x)(x-\mu_-)}}{2\pi x(1-x)}dx$\\
\end{tabular} & $I_W = [\mu_-,\mu_+]$ & \rule{0 in}{.5in}\begin{tabular}{c} $\mu_{\pm} = \left(\displaystyle\frac{\displaystyle\sqrt{b}\pm \displaystyle\sqrt{a(a+b-1)}}{a+b}\right)^2$, \\ $a,b\geq 1$\\
\end{tabular}\\ && \\ \hline
\end{tabular}\end{center}
\caption{Random matrix laws in raw form. The Kesten-McKay and Wachter laws are related by the linear transform $(2x_{\mbox {Wachter}}-1)v = x_{\mbox{Kesten-McKay}}$ and $a = b = v/2$. \label{tab:rawform}}
\end{table}

In raw form, these laws (Table \ref{tab:rawform})  appear as somewhat complicated expressions involving square roots. This paper highlights a unifying principle that underlies these four laws, namely the laws may be encoded as Jacobi symmetric tridiagonal matrices that are Toeplitz with a length 1 boundary.

This suggests that some of the nice properties of the big laws are connected to this property, and further suggests the importance of the larger family of laws encoded as Toeplitz with length $k$ boundary, known as  ``nearly Toeplitz" matrices.  This motivates the two parts of this paper:

\begin{enumerate}
\item  We tabulate in one place key properties of the four laws, not all of which can be found in the literature. These sections are expository, with the exception of the as-of-yet unpublished Wachter moments, and the Kesten-McKay and Wachter law Jacobi parameters and free cumulants.
\item We describe a new algorithm to exploit the Toeplitz-with-length-$k$ boundary structure.
In particular, we show how practical it is  to approximate distributions with incomplete information using distributions having nearly-Toeplitz encodings.
 \end{enumerate} 

Studies of nearly Toeplitz matrices in random matrix theory have been pioneered by Anshelevich \cite{Anshelevich2011,Anshelevich2011b}.

Historically, the Wigner semicircle law is the most famous.  The weight function is classical, and corresponds to Chebychev polynomials of the second kind. It is the equilibrium measure \cite{Deift1998} for Hermite Polynomials and the asymptotic distribution for Gaussian or Hermite ensembles (GOE, GUE, GSE, etc.).  None of the other weight functions are classical, but they are all equilibrium measures for classical polynomials. The second most famous law is the Marchenko-Pastur law. It is the equilibrium measure for Laguerre Polynomials and is the asymptotic distribution for Wishart matrices or Laguerre ensembles.  The Kesten-McKay law, described in \cite{Hora2007}, is the equilibrium measure for Gegenbauer Polynomials. It is not commonly included among the Wigner semicircle, Marchenko-Pastur, and Wachter laws, but we believe that it merits inclusion on account of its place in the upper-right box in Table \ref{tab:jacobi}. The Wigner and arcsine distributions are special cases of the Kesten-McKay distribution. The Wachter Law generalizes to Jacobi Polynomials. They describe sections of random unitary matrices, MANOVA matrices, and general Jacobi ensembles.

The laws have been encoded in many formats over the years. Wigner's earliest work encoded the semicircle law through its exponential generating function, the Bessel function. For many decades the Stieltjes or Cauchy transforms (Table \ref{tab:cauchy})  have been valuable.  Free probability has proposed the R-transform and S-transform as efficient encodings for these laws (Table \ref{tab:RS}).  In this note, we list each representation of each distribution, but we focus on the Jacobi parameter representation.
 
Other laws may be characterized as being asymptotically Toeplitz or numerically Toeplitz fairly quickly, such as the limiting histogram of the eigenvalues of $(X/\sqrt{m}+\mu I)^t(X/\sqrt{m}+\mu I)$, where $X$ is $m\times n$, $n$ is $O(m)$, and $m\longrightarrow\infty$ (Figure \ref{fig:cubic}).

This nearly Toeplitz property inspires an underlying approximation concept. Instead of simply truncating a Jacobi matrix, we can construct an infinite Toeplitz matrix that is Toeplitz on all but an initial finite set of rows and columns. In turn, we apply this approximation idea to the moment problem in Section \ref{sec:moment}: given a finite set of moments, we can use the Lanczos iteration to find corresponding Jacobi parameters, put those parameters in a continued fraction, and take its inverse Stieltjes transform to find a smooth, compactly-supported distribution with the correct moments. Instead of using the moments to find the Jacobi parameters we can also use a discretization of the measure. In \cite{Rao2007}, the Cauchy transform and its relationship to continued fractions of depth $2$ are discussed, an idea that is generalized here in Table \ref{tab:algorithm}.

Mathematical investigations into weight functions with Toeplitz or asymptotically Toeplitz Jacobi parameters may be found in the work of two mathematicians with coincidentally similar names: Geronimus \cite{Geronimus1977} and Geronimo \cite{Geronimo1980}, \cite{Geronimo2014}.  It seems known that the algorithm may have issues, perhaps reminiscent of the Gibbs phenomenon of Fourier analysis, in that atoms may emerge.

In retrospect, a good part of this paper may be found explicitly or implicitly in the work of Anshelevich \cite{Anshelevich2011b}.  Nonetheless, as we began to form the various tables  and noticed the nice Catalan and Narayana properties, and especially the little box in the upper right in Table \ref{tab:jacobi}, we realized that in addition to the algorithm, there were enough ideas that were below the surface that we worked out for ourselves and wished to share.

\section{The Jacobi Symmetric Tridiagonal Encoding of Probability Distributions}

All distributions have corresponding tridiagonal matrices of Jacobi parameters. They may be computed, for example, by the continuous Lanczos iteration, described in \cite[p.286]{Trefethen1997} and reproduced in Table \ref{tab:lanczos}.

We computed  the Jacobi representations of the four laws providing the results in Table \ref{tab:jacobi}. The Jacobi parameters ($\alpha_i$ and $\beta_i$ for $i = 0,1,2,\ldots$) are elements of an  infinite Toeplitz tridiagonal representations bordered by the first row and column, which may have different values from the Toeplitz part of the matrix.

\begin{table}[h!]
\rule{.1in}{0in}
$
\left[\begin{array}{c|ccccc}
\alpha_0 & \beta_0 & & & & \\
\hline
\beta_0 & \alpha_1 & \beta_1 & & & \\
 & \beta_1 & \alpha_1 & \beta_1 & & \\
& & \ddots & \ddots & \ddots & \\
& & & \beta_1 & \alpha_1 & \beta_1 \\
& & & & \beta_1 & \alpha_1
\end{array}\right]
$
%
\rule{.1in}{0in}\begin{tabular}{|c|c|c|}
\hline
 & $\alpha_0 = \alpha_1$ & $\alpha_0 \neq \alpha_1$ \\
\hline
$\beta_0 = \beta_1$ & Wigner semicircle & Marchenko-Pastur \\
\hline
$\beta_0\neq\beta_1$ & Kesten-McKay & Wachter \\
\hline
\end{tabular}
\begin{center}
\begin{tabular}{|lcccc|}
\hline
\rule{-.1 in}{.25in}
Measure & $\alpha_0$ & $\alpha_n$, $(n\geq 1)$ & $\beta_0$ & $\beta_n$, $(n\geq 1)$\\
&&&&\\
\hline
\rule{-.1 in}{.35in}
Wigner Semicircle & $0$ & $0$ & $1$ & $1$\\
&&&&\\
\hline
\rule{-.1 in}{.35in}
Marchenko-Pastur & $\lambda$  & $\lambda+1$ & $\sqrt{\lambda}$  & $\sqrt{\lambda}$ \\
&&&&\\
\hline
\rule{-.1 in}{.35in}
Kesten-McKay & $0$ & $0$ & $\sqrt{v}$ & $\sqrt{v-1}$ \\
&&&&\\
\hline
\rule{-.1 in}{.35in}
Wachter & $\displaystyle\frac{a}{a+b}$ & $\displaystyle\frac{a^2-a+ab+b}{(a+b)^2}$ & $\displaystyle\frac{\sqrt{ab}}{(a+b)^{3/2}}$ & $\displaystyle\frac{\sqrt{ab(a+b-1)}}{(a+b)^2}$  \\
&&&&\\ \hline
\end{tabular}
\end{center}
\caption{Jacobi parameter encodings for the big level density laws.  Upper left: Symmetric  Toeplitz  Tridiagonal with 1-boundary , Upper Right:  Laws Organized by Toeplitz Property, Below: Specific Parameter Values \label{tab:jacobi}}
\end{table}

Anshelovich \cite{Anshelevich2011b} provides a complete table of six distributions that have Toeplitz Jacobi structure. The first three of which are semicircle, Marchenko-Pastur, and Wachter.  The other three distributions occupy the same box as Wachter in Table \ref{tab:jacobi}. Anshelovich casts the problem as the description of all distributions whose orthogonal polynomials have generating functions of the form
$$ \sum_{n=0}^{\infty}P_n(x)z^n = \frac{1}{1-xu(z)+tv(z)}, $$
which he calls Free Meixner distributions.

He includes the one and two atom forms of the Marchenko-Pastur and Wachter laws which correspond in random matrix theory to the choices of tall-and-skinny vs. short-and-fat matrices in the SVD or CS decompositions, respectively.

\section{Infinite RMT Laws.}

This section compares the properties of all four major infinite random matrix theory laws, the Wigner semicircle law, the Marchenko-Pastur law, the Kesten-McKay law, and the Wachter law.

We state the four laws of infinite dimensional random matrix theory and their intervals of support. $\mu_{WS}$ on $I_{WS}$, $\mu_{MP}$ on $I_{MP}$, $\mu_{KM}$ on $I_{KM}$, and $\mu_W$ on $I_W$ correspond to the Wigner semicircle, Marchenko-Pastur, Kesten-McKay, and Wachter laws. These laws are originally credited  to \cite{Wigner1958}, \cite{Marcenko1967}, \cite{McKay1981}, and \cite{Wachter1978}. See Table \ref{tab:rawform}.

\begin{table}[h!]
\begin{center}
\begin{tabular}{|lc|}
\hline
\rule{-.1 in}{.25in}
Measure & Cauchy Transform\\
&\\
\hline
\rule{-.1 in}{.35in}
Wigner Semicircle & $\displaystyle\frac{z-\sqrt{z^2-4}}{2}$  \\
&\\
\hline
\rule{-.1 in}{.35in}
Marchenko-Pastur &  $\displaystyle\frac{1-\lambda+z-\sqrt{(1-\lambda+z)^2-4z}}{2z}$\\
&\\
\hline
\rule{-.1 in}{.35in}
Kesten-McKay & $\displaystyle\frac{(v-2)z-v\sqrt{4(1-v)+z^2}}{2(v^2-z^2)}$ \\
&\\
\hline
\rule{-.1 in}{.35in}
Wachter & $\displaystyle\frac{1-a+(a+b-2)z-\sqrt{(a+1-(a+b)z)^2-4a(1-z)}}{2z(1-z)}$ \\
&\\
\hline
\end{tabular}
\end{center}
\vspace{-.2in}
\caption{Cauchy transforms \label{tab:cauchy}.} \vspace{.2in}
\end{table}

We can also write down the moments for each measure in Table \ref{tab:moments}, for Wigner and Marchenko-Pastur see \cite{Edelman2014}, for Kesten-McKay see \cite{McKay1981}, and for Wachter see Theorem 6.1 in the Section 6. Remember the Catalan number $C_n = \frac{1}{n+1}\binom{2n}{n}$ and the Narayana polynomial $N_{n}(r)=\sum_{j=1}^{n}N_{n,j}r^{j},$ where $N_{n,j} = \frac{1}{n}\binom{n}{j}\binom{n}{j-1}$, excepting $N_0(r) = 1$. The coefficients of $v^j(1-v)^{n/2-j}$ in the Kesten-McKay moments form the Catalan triangle. We discuss the pyramid created by the Wachter moments in Section \ref{sec:moments}.

\begin{table}[h!]
\begin{center}
\rule{-.2 in}{0in}\begin{tabular}{|lc|}
\hline
\rule{-.1 in}{.25in} Measure & Moment $n$\\
&\\
\hline
\rule{-.1 in}{.35in} Wigner Semicircle & $C_{n/2}$ if $n$ is even, $0$ otherwise  \\
&\\
\hline
\rule{-.1 in}{.35in} Marchenko-Pastur & $N_n(\lambda)$ \\
&\\
\hline
\rule{-.1 in}{.45in} Kesten-McKay & $\displaystyle\sum_{j=1}^{n/2}\binom{n-j}{n/2}\left(\frac{j}{n-j}\right)v^j(v-1)^{n/2-j}$ if $n$ is even, $0$ otherwise\\
&\\
\hline
\rule{-.1 in}{.45in} Wachter &  $\displaystyle\frac{a}{a+b} - (a+b)\displaystyle\sum_{j=0}^{n-2}\left[\left(\frac{\sqrt{a(a+b-1)}}{a+b}\right)^{2j+4}N_{j+1}\left(\frac{b}{a(a+b-1)}\right)\right]$\\
&\\ \hline
\end{tabular}
\end{center}
\vspace{-.2in}
\caption{Moments \label{tab:moments}}
\end{table}

Inverting the Cauchy transforms and subtracting $1/w$,  computes the $R$-transform, see Table \ref{tab:RS}. If there are multiple roots,  we pick one with a series expansion with no pole at $w=0$.

\begin{table}[h!]
\begin{center}
\begin{tabular}{|lcc|}
\hline
\rule{-.1 in}{.25in}
Measure & $R$-transform & $S$-transform \\
&&\\
\hline
\rule{-.1 in}{.35in}
Wigner Semicircle &  $w$  & 1\\
&&\\
\hline
\rule{-.1 in}{.35in}
Marchenko-Pastur &  $\displaystyle\frac{\lambda}{1-w}$ &  $\displaystyle\frac{z-\lambda}{z^2}$ 
 \\
&&\\
\hline
\rule{-.1 in}{.35in}
Kesten-McKay & $\displaystyle\frac{-v + v\sqrt{1+4w^2}}{2w}$&  $\displaystyle\frac{v}{v^2-z^2}$\\
& & \\
\hline
\rule{-.1 in}{.35in}
Wachter & $\displaystyle\frac{-a-b+w+\sqrt{(a+b)^2+2(a-b)w + w^2}}{2w}$ & $\displaystyle\frac{a-az-bz}{z^2(z-1)}$ \\
&&\\ \hline
\end{tabular}
\end{center}
\caption{$R$-transforms  and  $S$-transforms  computed as $S(z) = R^{-1}(z)/z$ \label{tab:RS} .}
\vskip .2 in
\end{table}

The free cumulants $\kappa_n$ for each measure appear in Table \ref{tab:cumulants} by expanding the $R$-transform above (the generating function for the Narayana polynomials is given by \cite{Mansour2009}, the generating function for the Catalan numbers is well known).

It is widely known that the Catalan numbers are the moments of the semicircle law, but we have not seen any mention that the same numbers figure prominently as the free cumulants of the Kesten-McKay Law. The Narayana Polynomials are prominent as the moments of the Marchenko-Pastur Law, but they also figure clearly as the free cumulants of the Wachter Law. There are well known relationships, involving Catalan numbers, between the moments and free cumulants of any law \cite{Nica2006}, but we do not know if the pattern is general enough to take the moments of one law, transform it somewhat, and have them show up in the free cumulants in another law.

\begin{table}[h!]
\begin{center}
\begin{tabular}{|lc|}
\hline
\rule{-.1 in}{.25in}
Measure & $\kappa_n$\\
& \\
\hline
\rule{-.1 in}{.35in}
Wigner Semicircle & $\delta_{n,2}$ \\
& \\
\hline
\rule{-.1 in}{.35in}
Marchenko-Pastur & $\lambda$ \\
& \\
\hline
\rule{-.1 in}{.35in}
Kesten-McKay &  $(-1)^{(n-2)/2}vC_{(n-2)/2}$ if $n$ is even, $0$ otherwise\\
&\\
\hline
\rule{-.1 in}{.35in}
Wachter & $-N_n\left(-\displaystyle\frac{b}{a}\right)\displaystyle\frac{(-a)^{n+1}}{(a+b)^{2n+1}}$ \\
& \\ \hline
\end{tabular}
\end{center}
\caption{Free cumulants. \label{tab:cumulants}}
\vskip .2 in
\end{table}

We compute an S-transform as  $S(z) = R^{-1}(z)/z$. See Table \ref{tab:RS}.

Each measure has a corresponding three-term recurrence for its orthonormal polynomial basis, with $q_{-1}(x) = 0$, $q_0(x) = 1$, $\beta_{-1} = 0$, and for $n\geq 0$, $q_{n+1}(x) = ((x-\alpha_n)q_{n}(x)-\beta_{n-1}q_{n-1}(x))/\beta_n$. In the case of the Wigner semicircle, Marchenko-Pastur, Kesten-McKay, and Wachter laws, the Jacobi parameters $\alpha_n$ and $\beta_n$ are constant for $n\geq 1$ because they are all versions of the Meixner law \cite{Anshelevich2011b} (a linear transformation may be needed). The Wigner Semicircle case is given by simplifying the Meixner law in \cite{Anshelevich2011}, and the Marchenko-Pastur, Kesten-McKay, and Wachter cases are given by taking two iterations Lanczos algorithm symbolically to get $\alpha_1$ and $\beta_1$. See Table \ref{tab:jacobi}.

\begin{table}[h!]
\begin{center}
\rule{-.05in}{0in}\begin{tabular}{|lc|}
\hline
\rule{-.1 in}{.25in}
Measure & $q_n(x)$, $n\geq 1$.\\
&\\
\hline
\rule{-.1 in}{.35in}
Wigner Semicircle & $U_n\left(\frac{x}{2}\right)$ \\
&\\
\hline
\rule{-.1 in}{.35in}
Marchenko-Pastur &  $\lambda^{(n-1)/2}(x-\lambda)U_{n-1}\left(\frac{x-\lambda-1}{2\sqrt{\lambda}}\right) - \lambda^{n/2}U_{n-2}\left(\frac{x-\lambda-1}{2\sqrt{\lambda}}\right)$ \\
&\\
\hline
\rule{-.1 in}{.35in}
Kesten-McKay & $ (v-1)^{(n-1)/2}xU_{n-1}\left(\frac{x}{2\sqrt{v-1}}\right) - v(v-1)^{(n-2)/2}U_{n-2}\left(\frac{x}{2\sqrt{v-1}}\right)$\\
&\\
\hline
\rule{-.1 in}{.35in}
Wachter & $ \left(x-\frac{a}{a+b}\right)\left(\frac{\sqrt{ab(a+b-1)}}{(a+b)^2}\right)^{n-1}U_{n-1}\left(\frac{-b - a (a + b-1) + (a + b)^2 x}{2 \sqrt{a b ( a + b-1)}}\right)$ \\
& $- \frac{a+b}{a+b-1}\left(\frac{\sqrt{ab(a+b-1)}}{(a+b)^2}\right)^n U_{n-2}\left(\frac{-b - a (a + b-1) + (a + b)^2 x}{2 \sqrt{a b (a + b-1})}\right) $\\
&\\ \hline
\end{tabular}
\end{center}
\caption{Sequences of polynomials orthogonal over of the four major laws. \label{tab:orthogonal}}
\end{table}

Each measure also has an infinite sequence of monic polynomials $q_n(x)$ which are orthogonal with respect to that measure. They can be written as sums of Chebyshev polynomials of the second kind, $U_n(x)$, which satisfy $U_{-1} = 0$, $U_0(x) = 1$, and $U_n(x) = 2xU_{n-1}(x)-U_{n-2}(x)$ for $n\geq 1$, \cite{Kusalik2005}. See Table \ref{tab:orthogonal}. For $n = 0$, $q_0(x) = 1$, and in general for $n\geq 1$,
$$ q_n(x) = \beta_1^{n-1}(x-\alpha_0)U_{n-1}\left((x-\alpha_1)/(2\beta_1)\right) - \beta_0^2\beta_1^{n-2}U_{n-2}\left((x-\alpha_1)/(2\beta_1)\right).  $$
In the Wigner semicircle case the polynomials can be combined using the recursion rule for Chebyshev polynomials.

\section{The Wachter Law Moment  Pyramid.}
\label{sec:moments}

\begin{figure}
\centering
\includegraphics[scale=.6]{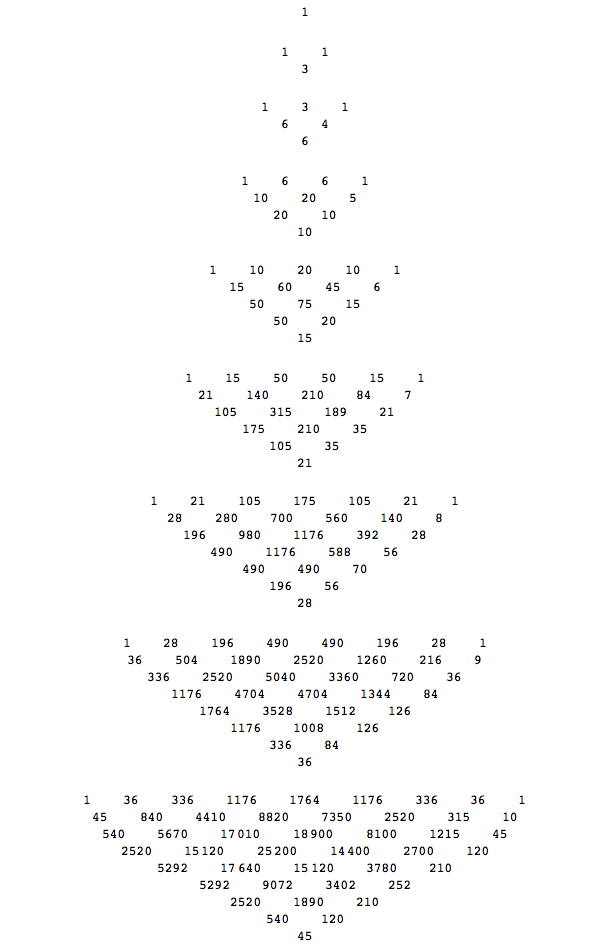}
\caption{A number pyramid from the coefficients of the  Wachter law  moments. \label{fig:wachtermoments}}
\end{figure}

Using Mathematica we can extract an interesting number pyramid from the Wachter moments, see Figure \ref{fig:wachtermoments}. Each triangle in the pyramid is formed by taking the coefficients of $a$ and $b$ in the $i$-th Wachter moment, with the row number within the pyramid determined by the degree of the corresponding monomial in $a$ and $b$. All factors of $(a+b)$ are removed from the numerator and denominator beforehand and alternating signs are ignored.

Furtheremore, there are many patterns within the pyramid. The top row of each triangle is a list of Narayana numbers, which sum to Catalan numbers. The bottom entries of each pyramid are triangular numbers. The second-to-bottom entry on the right of every pyramid is a sum of consecutive triangular numbers. The second to both the left and right on the top row of every triangle are also triangular numbers.

\section{Moments Build Nearly-Toeplitz Jacobi Matrices.}
\label{sec:moment}

This section is concerned with recovering a probability distribution from its Jacobi parameters, $\alpha_i$ and $\beta_i$ such that they are ``Nearly Toeplitz,'' i.e. there exists a $k$ such that for $i\geq k$ all $\alpha_i$ are equal and all $\beta_i$ are equal. Note that $i$ ranges from $k$ to $\infty$. The Jacobi parameters are found from a distribution by the Lanczos iteration.

We now state the continuous Lanczos iteration, replacing the matrix $A$ by the variable $x$ and using $\mu = \mu_{WS}, \mu_{MP}, \mu_M, \mu_W$ to compute dot products. A good source is \cite{Trefethen1997}. For a given measure $\mu$ on an interval $I$, let
$$ (p(x),q(x)) = \int_I p(x)q(x)d\mu, $$
and $\|p(x)\| = \sqrt{(p(x),p(x))}$. Then the Lanczos iteration is described by Table \ref{tab:lanczos}.

\begin{table}[h!]
\begin{center}
\fbox{
\begin{minipage}{4.4 in}
\begin{center}
\underline{Lanczos on Measure $\mu$}
\end{center}
\begin{algorithmic}
\State $\beta_{-1} = 0$, $q_{-1}(x) = 0$, $q_0(x) = 1$
\For{$n = 0,1,2,\ldots}$
	\State $v(x) = xq_n(x)$
	\State $\alpha_n = (q_n(x),v(x))$
	\State $v(x) = v(x) - \beta_{n-1}q_{n-1}(x)-\alpha_nq_n(x)$
	\State $\beta_n = \|v(x)\|$
	\State $q_{n+1}(x) = v(x)/\beta_n$
\EndFor
\end{algorithmic}
\end{minipage}
} \end{center}
\caption{The Lanczos iteration produces the Jacobi parameters in $\alpha$ and $\beta$. \label{tab:lanczos}}
\end{table}

There are two ways to compute the integrals numerically. The first is to sample $x$ and $q_n(x)$ at many points on the interval of support for $q_0(x) = 1$ and discretize the integrals on that grid. The second can be done if you know the moments of $\mu$. If $r(x)$ and $s(x)$ are polynomials, $(r(x),s(x))$ can be computed given $\mu$'s moments. Since the $q_n(x)$ are polynomials, every integral in the Lanczos iteration can be done in this way. In that case, the $q_n(x)$ are stored by their coefficients of powers of $x$ instead of on a grid. Once we have reached $k$ iterations, we have fully constructed the infinite Jacobi matrix using the first batch of $\mu$'s moments, or a discretization of $\mu$.

\begin{table}[h!]
{\bf Algorithm: Compute Measure from Nearly Toeplitz Jacobi Matrix.}
\begin{enumerate}
\item Nearly Jacobi Toeplitz Representation: Run the continuous Lanczos algorithm up to step $k$, after which all $\alpha_i$ are equal and all $\beta_i$ are equal, or very nearly so. If they are equal, this algorithm will recover $d\mu$ exactly, otherwise it will find it approximately. The Lanczos algorithm may be run using a discretization of the measure $\mu$, or its initial moments.
$$ (\alpha_{0:\infty},\beta_{0:\infty}) = {\rm Lanczos}\left(d\mu(x)\right). $$
\item  Cauchy transform: evaluate the finite continued fraction below on the interval of $x$ where it is imaginary.
$$ g(x) = \cfrac{1}{x-\alpha_0-\cfrac{\beta_0^2}{x-\alpha_1-\cfrac{\beta_1^2}{\ddots -\cfrac{\beta_{k-2}^2}{\alpha_{k-1}-\cfrac{2\beta_{k-1}^2}{x-\alpha_k+\sqrt{(\alpha_k-x)^2-4\beta_k^2}}}}}}. $$
\item Inverse Cauchy Transform:   divide the imaginary part by $-\pi$, to compute the desired measure.
$$ d\mu(x) = -\frac{1}{\pi}{\rm Im}\left(g(x)\right). $$
\end{enumerate}
\caption{Algorithm recovering or approximating an analytic measure by Toeplitz matrices with boundary. \label{tab:algorithm}}
\end{table}

Step 1 can start with a general measure in which case Step 3 finds an approximate measure with a nearly Toeplitz representation. Step 1 could also start with a sequence of moments. It should be noted that the standard way to go from moments to Lanczos coefficients uses a Hankel matrix of moments and its Cholesky factorization (\cite{Golub1969}, (4.3)).

\begin{figure}
\centering
\includegraphics[scale=.8]{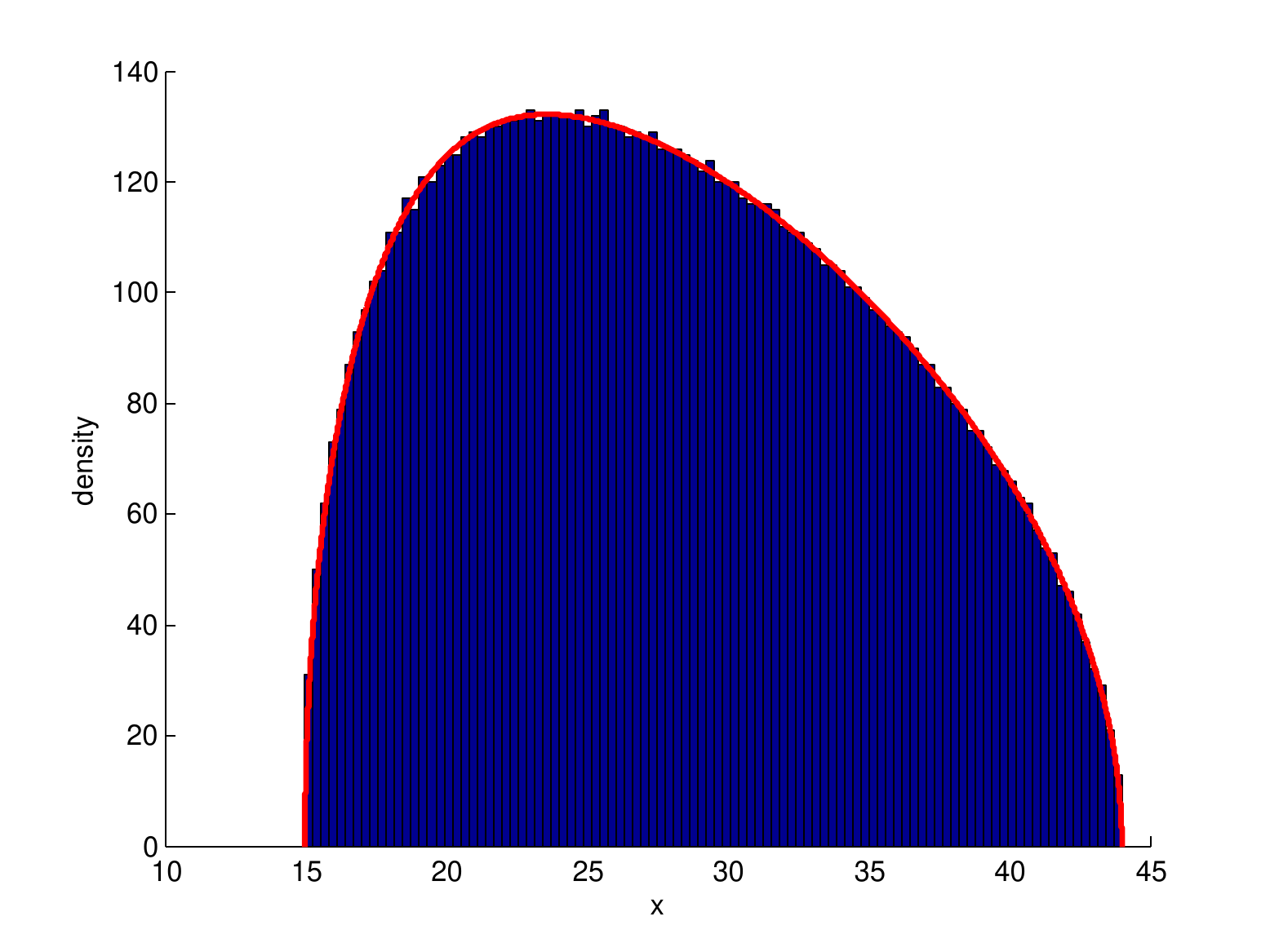}
\caption{Eigenvalues taken from $(X/\sqrt{m}+\mu I)^t(X/\sqrt{m}+\mu I)$, where $X$ is $m\times n$, $m = 10^4$, $n=3m$, $\mu = 5$. The blue bar histogram is taken using hist.m, a better one was taken by convolving the data with Gaussian kernel. That convolution histogram was used to initialize the continuous Lanczos algorithm which produced five $\alpha$'s and $\beta$'s. They were put into a continued fraction as described above, assuming $\alpha_i$ and $\beta_i$ to be constant after $i = 5$. The continued fraction recreated the histogram, which is the thick red line.
\label{fig:cubic}}
\end{figure}

\begin{figure}
\centering
\includegraphics[scale=.6]{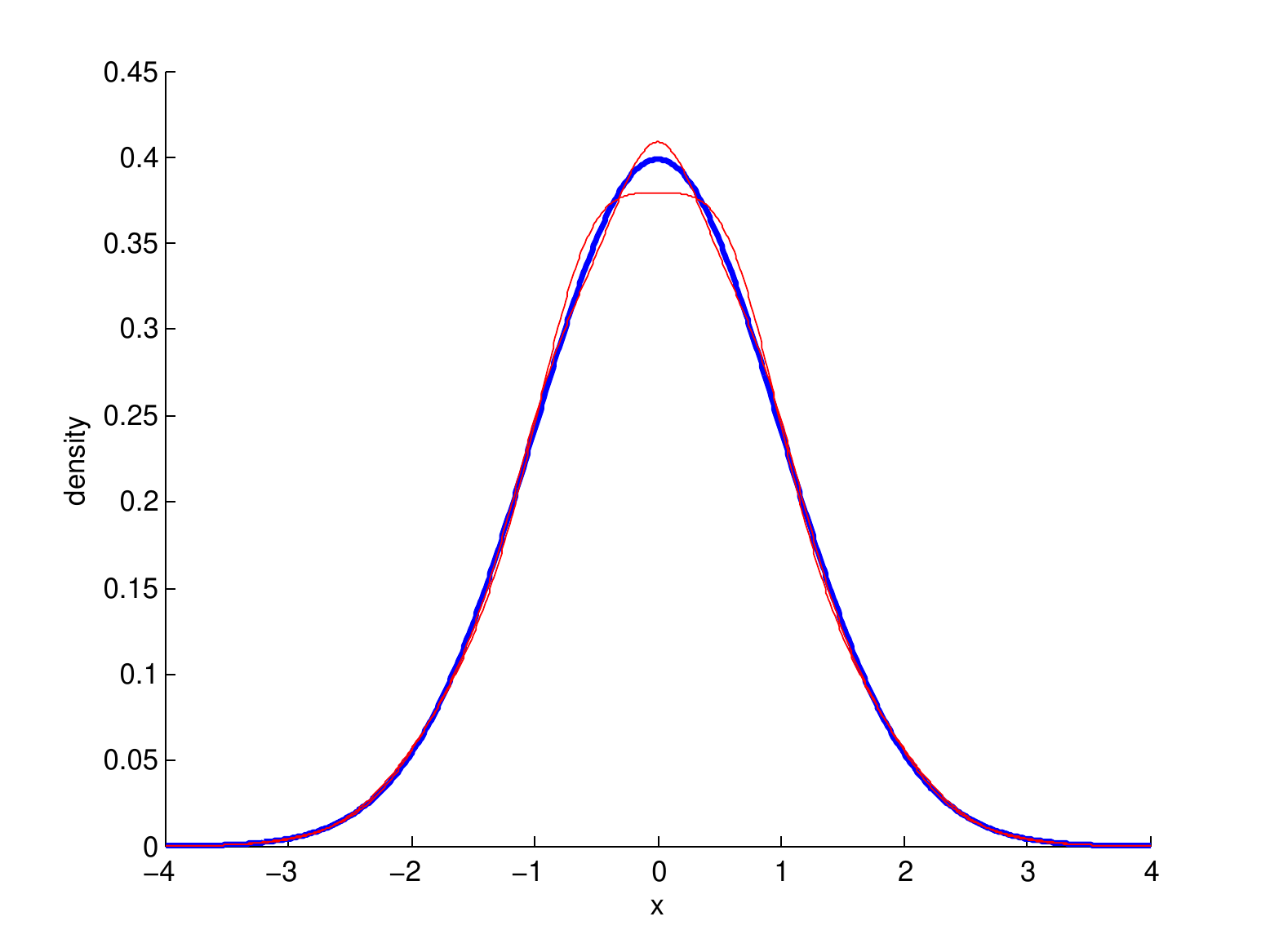}
\caption{The normal distribution's Jacobi matrix is not well approximated by Toeplitz plus boundary, but with sufficiently many moments good approximations are possible. The above graph shows the normal distribution recovered by the method in this paper using $10$ and $20$ moments. The thick line is the normal computed by $e^{-x^2/2}\sqrt{2\pi}$, and the thin lines on top of it use our algorithm.  \label{fig:normal}
}
\end{figure}

As an example, we apply the algorithm to the histogram of the eigenvalues of $(X/\sqrt{m}+\mu I)^t(X/\sqrt{m}+\mu I)$, where $X$ is $m\times n$, which has Jacobi parameters $\alpha_i$ and $\beta_i$ that converge asymptotically and quickly. We smooth the histogram using a Gaussian kernel and then compute its Jacobi parameters. The reconstruction of the histogram is in Figure \ref{fig:cubic} We also use the above algorithm to reconstruct a normal distribution from its first sixty moments, see Figure \ref{fig:normal}.

The following theorem concerning continued fractions allows one to stably recover a distribution from its Lanczos coefficients $\alpha_i$ and $\beta_i$. As we have said, if the first batch of $\mu$'s moments are known, we can find all $\alpha_i$ and $\beta_i$ from $i=0$ to $\infty$ using the continuous Lanczos iteration.

\begin{theorem} Let $\mu$ be a measure on interval $I \subset \mathbb{R}$ with Lanczos coefficients  $\alpha_i$ and $\beta_i$, with the property that all $\alpha_i$ are equal for $i\geq k$ and all $\beta_i$ are equal for $i\geq k$. We can recover $I = [\alpha_k-2\beta_k,\alpha_k+2\beta_k]$, and we can recover $d\mu(x)$ using a continued fraction. This theorem combines Theorems 1.97 and 1.102 of \cite{Hora2007}.
$$ g(x) = \cfrac{1}{x-\alpha_0-\cfrac{\beta_0^2}{x-\alpha_1-\cfrac{\beta_1^2}{\ddots -\cfrac{2\beta_{k-1}^2}{x-\alpha_k+\sqrt{(\alpha_k-x)^2-4\beta_k^2}}}}} $$
$$ d\mu(x) = -\frac{1}{\pi}{\rm Im}\left(g(x)\right). $$
\end{theorem}

\begin{figure}
\centering
\includegraphics[scale=.6]{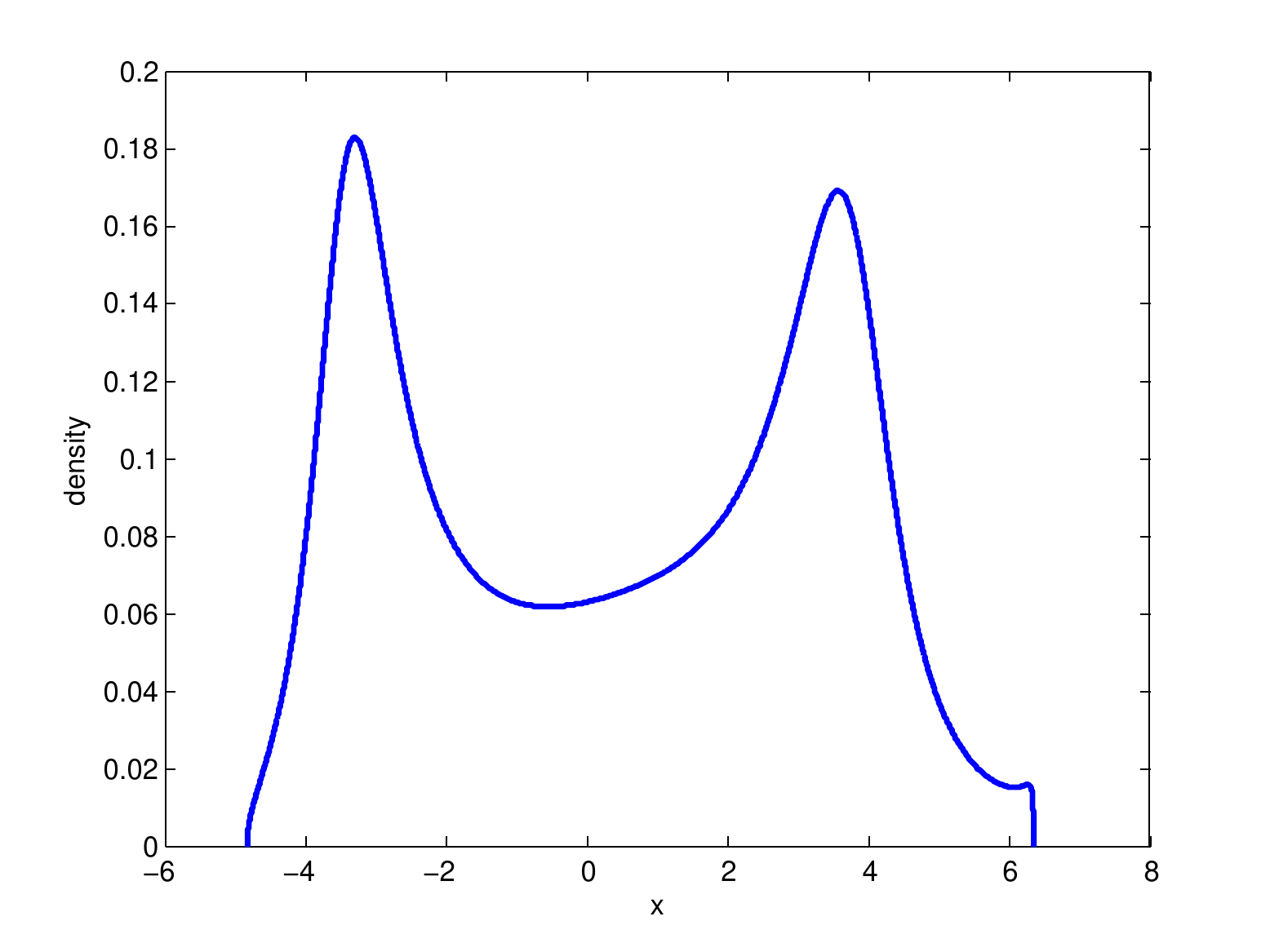}
\includegraphics[scale=.6]{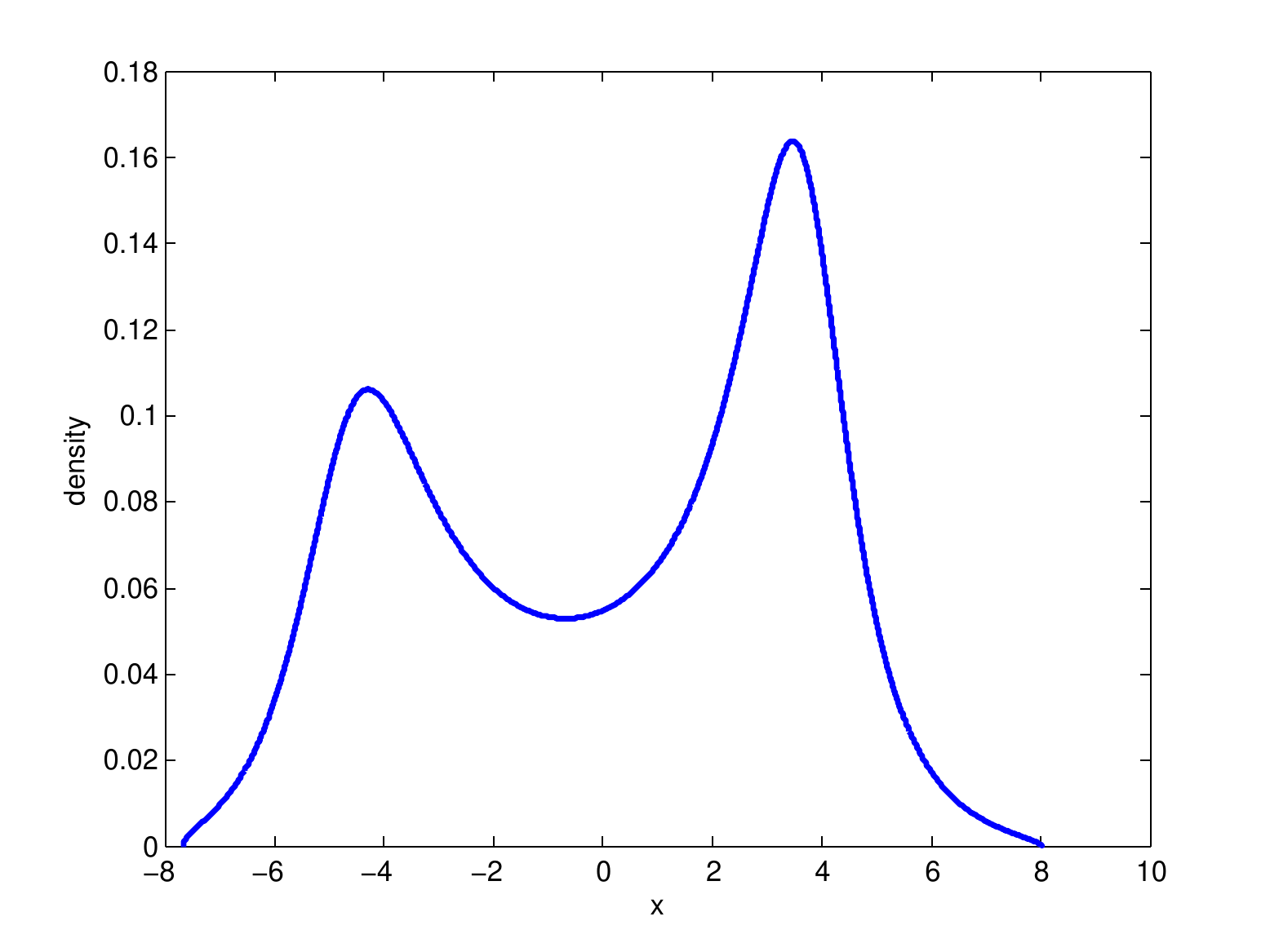}
\caption{Recovery from of a distribution from random $\alpha_i$ and $\beta_i$ using Theorem 5.1.  On top we use $k=5$, on bottom we use $k = 3$. \label{fig:recovery}}
\end{figure}

Figure \ref{fig:recovery} illustrates curves recovered from random terminating continued fractions $g(x)$ such that the $\beta_i$ are positive and greater in magnitude than the $\alpha_i$. In both cases, the above theorem allows correct recovery of the $\alpha_i$ and $\beta_i$ (which is not always numerically possible). In the first one, $k = 5$, in the second, $k = 3$.

If $X$ is an $m\times n$, $m < n$ matrix of normals for $m$ and $n$ very large, $(X/\sqrt{m}+\mu I)^t(X/\sqrt{m}+\mu I)$ has $\alpha_i$ and $\beta_i$ which converge to a constant, making its eigenvalue distribution recoverable up to a very small approximation. See Figure \ref{fig:cubic}

We also tried to reconstruct the normal distribution, whose Jacobi parameterization is not at all Toeplitz, and which is not compactly supported. Figure \ref{fig:normal} plots the approximations using $10$ and $20$ moments.

\section{Direct computation of the Wachter law moments.}

While the moments of the Wachter law may be obtained in a number of ways, including expanding the Cauchy Transform, or applying the mobius inverse formula to the free cumulants, in this section we show that a direct computation of the integral is possible.

\begin{theorem} We find the moments of the Wachter law, $m_k$.
$$ m_k = \frac{a}{a+b} - (a+b)\sum_{j=0}^{k-2}\left[\left(\frac{\sqrt{a(a+b-1)}}{a+b}\right)^{2j+4}N_{j+1}\left(\frac{b}{a(a+b-1)}\right)\right].$$
\end{theorem}
\begin{proof}
We start by integrating the following expression by comparing it to the Marchenko-Pastur law.
$$ J_1 = \frac{1}{2\pi}\int_{\mu_-}^{\mu_+}x^k\sqrt{(\mu_+-x)(x-\mu_-)}dx. $$
If $x = su$, $dx = sdu$ and this integral becomes
$$ J_1 = \frac{s^{k+2}}{2\pi}\int_{\frac{\mu_-}{s}}^{\frac{\mu_+}{s}}u^k\sqrt{\left(\frac{\mu_+}{s}-u\right)\left(u-\frac{\mu_-}{s}\right)}du. $$
To compare this expression to the Marchenko-Pastur law, we need to pick $s$ and $\lambda$ such that $\frac{\mu_+}{s} = (1+\sqrt{\lambda})^2$ and $\frac{\mu_-}{s} = (1-\sqrt{\lambda})^2$ for $\lambda \geq 1$. There are more than one choices of each parameter, but we pick $\sqrt{s} = \frac{1}{2}\left(\sqrt{\mu_+}-\sqrt{\mu_-}\right)$ and
$$ \sqrt{\lambda} = \frac{\sqrt{\mu_+}+\sqrt{\mu_-}}{\sqrt{\mu_+}-\sqrt{\mu_-}}. $$
Using the Narayana numbers, and the formula for the moments of the Marchenko-Pastur law, the integral equals
$$ J_1 = \left(\frac{1}{2}\left(\sqrt{\mu_+}-\sqrt{\mu_-}\right)\right)^{2k+4}N_{k+1}\left(\left(\frac{\sqrt{\mu_+}+\sqrt{\mu_-}}{\sqrt{\mu_+}-\sqrt{\mu_-}}\right)^2\right). $$
Using $a$ and $b$, this becomes
$$ J_1 = \left(\frac{\sqrt{a(a+b-1)}}{a+b}\right)^{2k+4}N_{k+1}\left(\frac{b}{a(a+b-1)}\right). $$
We also need to integrate
$$ J_2 = \frac{1}{2\pi}\int_{\mu_-}^{\mu_+}\frac{\sqrt{(\mu_+-x)(x-\mu_-)}}{1-x}dx $$
Let $su = x-1$. $sdu = dx$. This becomes
$$J_2=-\frac{s}{2\pi}\int_{\frac{\mu_--1}{s}}^{\frac{\mu_+-1}{s}}\frac{\sqrt{\left(\frac{\mu_+-1}{s}-u\right)\left(u-\frac{\mu_--1}{s}\right)}}{u}du, $$
which by symmetry is
$$J_2=\frac{s}{2\pi}\int_{\frac{1-\mu_+}{s}}^{\frac{1-\mu_-}{s}}\frac{\sqrt{\left(\frac{1-\mu_-}{s}-u\right)\left(u-\frac{1-\mu_+}{s}\right)}}{u}du. $$
Using the same technique as previously, $\sqrt{s} = \frac{1}{2}\left(\sqrt{1-\mu_-}-\sqrt{1-\mu_+}\right)$ and
$$ \sqrt{\lambda} = \frac{\sqrt{1-\mu_+}+\sqrt{1-\mu_-}}{\sqrt{1-\mu_-}-\sqrt{1-\mu_+}}. $$
Using the fact that the Marchenko-Pastur law is normalized, the answer is
$$ J_2 = s = \frac{1}{4}\left(\sqrt{1-\mu_-}-\sqrt{1-\mu_+}\right)^2 = \frac{a}{(a+b)^2}.$$

Now we are ready to find the moments of the Wachter law. Using the geometric series formula,
\begin{eqnarray*}
m_k &=& \frac{a+b}{2\pi}\int_{\mu_-}^{\mu_+}\frac{x^{k-1}\sqrt{(\mu_+-x)(x-\mu_-)}}{1-x}dx \\
&=& \frac{a+b}{2\pi}\sum_{j=k-1}^{\infty}\int_{\mu_-}^{\mu_+}x^j\sqrt{(\mu_+-x)(x-\mu_-)}dx \\
&=& \frac{a+b}{2\pi}\int_{\mu_-}^{\mu_+}\frac{\sqrt{(\mu_+-x)(x-\mu_-)}}{1-x}dx - \frac{a+b}{2\pi}\sum_{j=0}^{k-2}\int_{\mu_-}^{\mu_+}x^j\sqrt{(\mu_+-x)(x-\mu_-)}dx \\
&=& \frac{a}{a+b} - (a+b)\sum_{j=0}^{k-2}\left[\left(\frac{\sqrt{a(a+b-1)}}{a+b}\right)^{2j+4}N_{j+1}\left(\frac{b}{a(a+b-1)}\right)\right].
\end{eqnarray*}
\end{proof}

\section{Acknowledgements}

We would like to thank Michael LaCroix, Plamen Koev,  Sheehan Olver and Bernie Wang for interesting discussions. We gratefully acknowledge the support of the National Science Foundation:  DMS-1312831, DMS-1016125, DMS-1016086.

\end{document}